\documentclass[12pt]{article}
\usepackage{a4wide,amsfonts,amssymb,theorem,amstext}

\makeatletter
\@addtoreset{equation}{section}
\@addtoreset{table}{section}

\makeatother

\theoremstyle{break}
\theorembodyfont{\rmfamily}
\newtheorem{theorem}{Theorem}[section]
\newtheorem{corollary}{Corollary}[section]
\newcommand{\ie}{\emph{i.e.}}

\begin{document}

\begin{center}
\begin{flushright}
April 17, 1998
\\
math.QA/9804084
\end{flushright}
\vspace*{0.5cm}
\begin{Large}
{\bf On the bicrossproduct structures for the
${\cal U}_\lambda(iso_{\omega_2\dots \omega_N}(N))$ family of algebras}
\\[0.5cm]
\end{Large}
J. C. P\'erez Bueno\footnote{E-mail: pbueno@lie.ific.uv.es}
\\[0.5cm]
{\sl Departamento de F\'{\i}sica Te\'orica and IFIC
\\
Centro Mixto Universidad de Valencia--CSIC
\\
E--46100 Burjassot (Valencia), Spain}
\end{center}

\begin{abstract}
It is shown that the family of deformed algebras
${\cal U}_\lambda(iso_{\omega_2\dots \omega_N}(N))$ has a different
bicrossproduct structure for each $\omega_a=0$ in analogy to the undeformed
case.
\end{abstract}

\section{Introduction}

Deformed algebras (usually called `quantum groups') have received great
attention since the original works of Drinfel'd, Jimbo and Faddeev,
Reshetikhin
and Takhtajan \cite{Dri:87,Jim:85,Jim:86,Fad.Res.Tak:89} which gave a
(unique) deformation procedure for simple Lie algebras.
However, the deformation of non-simple Lie algebras has been characterized by
the lack of a definite prescription and this explains why inhomogeneous
algebras do not have a unique deformation.

A possible approach to deforming non-simple algebras is by extending
the contraction of Lie algebras to the framework of deformed Hopf algebras,
an idea originally introduced by  Celeghini \emph{et al.}
\cite{Cel.Gia.Sor.Tar:91,Cel.Gia.Sor.Tar:92}.
As is well known, the standard \.In\"on\"u-Wigner \cite{Ino.Wig:53}
contraction of (simple) Lie algebras leads to
non-simple algebras which have a semidirect structure, where the ideal
is the abelianized part of the original algebra.
By introducing higher powers in the contraction parameter or, equivalently, by
performing two (or more) successive contractions it is also possible to arrive 
to algebras with a central extension structure.

This simple mechanism becomes difficult to implement for deformed algebras for 
which it is usually necessary to redefine the deformation parameter in terms of
the contraction one to have a well-defined contraction limit
\cite{Cel.Gia.Sor.Tar:91,Cel.Gia.Sor.Tar:92}.
This is the case, for instance, of the $\kappa$-Poincar\'e algebra
\cite{Luk.Now.Rue.Tol:91,Luk.Rue.Tol:94}
in which the deformation parameter $\kappa$ appears as a redefinition of the
original (adimensional) parameter $q$ of $so_q(3,2)$ in terms of the De Sitter
radius $R$.

A way to skip some of the problems of the standard contraction procedure 
for deformed algebras is to use the method of `graded' contractions.
This mechanism was put forward by Moody, Montigny and Patera
\cite{Mon.Pat:91,Moo.Pat:91} for Lie algebras and has been applied recently to
describe a large set of deformed Hopf algebras
\cite{Bal.Her.Olm.San:94,Bal.Her.Olm.San:95}.
The scheme provides the deformation of all motion algebras of flat affine
spaces in $N$ dimensions (the deformed Cayley-Klein (CK) algebras 
${\cal U}_\lambda(iso_{\omega_2\dots\omega_N}(N))$~\footnote{
The orthogonal Cayley-Klein family of
algebras are the Lie algebras of the motion groups of real spaces with a
projective metric \cite{Som:10,Yag.Roz.Yas:66}.})
including, the $\kappa$-Poincar\'e algebra in arbitrary dimensions, other
deformations of the Poincar\'e $N$-dimensional algebra, the Galilei algebra,
etc.

A different point of view to study inhomogeneous deformed algebras is
provided by Majid's \emph{bicrossproduct structure}
\cite{Maj:90,Maj:90b,Maj:95}.
In this construction we find the analogue of the Lie algebra semidirect
structure (and of the central extension structure in the more general case)
for Hopf algebras and provides, for this reason, an appropriate setting for the
study of deformations of inhomogeneous algebras.
This structure covers most of the deformed algebras obtained by
contraction but not all (see \cite{Azc.Bue:96}).
Thus, in the case of deformed algebras,
the correspondence between contraction and semidirect structure that
exists in the Lie algebra setting is not straightforward.

Nevertheless, the study of the particular algebras for which the structure of
bicrossproduct is present, turns out to be useful to understand its properties
because the deformation is mainly encoded in the
action and coaction mappings that characterize the bicrossproduct, whereas the
(two) Hopf algebras from which the bicrossproduct deformed algebra is
constructed are usually undeformed.
In the appropriate limit of the deformation parameter we obtain the
undeformed algebra, the coaction mapping is trivialized and the action mapping
is given by the Lie algebra commutators so that we recover the semidirect
product structure.
A particular example is the $\kappa$-Poincar\'e algebra
\cite{Luk.Now.Rue.Tol:91} the
bicrossproduct structure of which was found by Majid and Ruegg
\cite{Maj.Rue:94}.

Recently \cite{Azc.Olm.Bue.San:97} (see also \cite{Per:97})
has been shown that the whole family of
deformed inhomogeneous CK algebras 
${\cal U}_\lambda(iso_{\omega_2\dots \omega_N}(N))$
has a bicrossproduct structure, in analogy to the semidirect one that
appears after the contraction which goes from $so(p,q)$ to
$iso(p,q)$~\footnote{
Note that associated to each contraction there exists, in
the undeformed level, a semidirect product structure and, in the deformed
level, a possible bicrossproduct structure associated to it.
In this particular case the contraction $so(p,q)\to iso(p,q)$ gives rise to
a semidirect structure in which we have $p+q$ abelian (momentum) generators
and a (pseudo-)orthogonal group acting on them.}
and that it remains under all the possible graded contractions.
However, the question that naturally arises is whether these
contractions carry new bicrossproduct structures related to the semidirect
ones of the undeformed algebra which are the result of each contraction (see
(\ref{semidirect}) below).
We prove in this paper that this is indeed the case so that, for every
graded contraction in the inhomogeneous deformed CK family
${\cal U}_\lambda(iso_{\omega_2\dots \omega_N}(N))$,
we have an associated bicrossproduct structure.
The (a priori non-obvious) fact that all the possible semidirect product 
structures of the undeformed inhomogeneous CK algebras have a direct 
counterpart in the deformed case is the main result of this paper.

The paper is organized as follows.
In sec. \ref{sec2} we provide an account of the (undeformed) CK algebras and
their graded contractions.
In sec. \ref{sec3} some of the results in \cite{Azc.Olm.Bue.San:97} are
summarized.
They will permit us in sec. \ref{sec4} to show that for each possible
graded contraction in the CK family a new bicrossproduct structure arises.
Our results are illustrated at the end with some examples.

\section{Cayley-Klein algebras}
\label{sec2}

Let us start by recalling the definition of the orthogonal Cayley-Klein family
of algebras.
The (orthogonal) real Lie algebra $so(N+1)$ can be endowed with a
${\mathbb Z}_2^{\otimes N}$ grading group and corresponding to its graded
contractions we may
introduce a set of Lie algebras depending on $2^N-1$ real parameters
\cite{Her.San:96b}.
This set includes the original $so(N+1)$ algebra, all the possible
pseudo-orthogonal ones and many contracted algebras, as well as
the $N(N+1)/2$ dimensional abelian one.
The simplicity of the original $so(N+1)$ algebra is lost for \emph{arbitrary}
contractions
and different algebras in this set may have different properties
(as the number of independent Casimir operators).

However there exists a subfamily, the members of which share many
properties with the (parent) simple Lie algebra and hence may be called
`quasi-simple'.
This family, denoted by $so_{\omega_1\dots \omega_N}(N+1)$, is a set of
algebras characterized by $N$ real parameters $(\omega_1,\dots,\omega_N)$ and
corresponds to a natural subset of all possible graded contractions that may be
obtained from $so(N+1)$ (within this family we find, for instance, the
original $so(N+1)$ algebra, the $N$-dimensional Poincar\'e algebra, the
Euclidean algebra, etc.).
These algebras correspond exactly to the motion algebras of the
geometries of a real space with a projective metric in the Cayley--Klein sense
\cite{Som:10,Yag.Roz.Yas:66} and are therefore called CK orthogonal algebras.
Their non-zero brackets are
\begin{equation}
[{\mathbb J}_{ab},{\mathbb J}_{ac}] =  \omega_{ab} {\mathbb J}_{bc}\quad ,
\quad
[{\mathbb J}_{ab},{\mathbb J}_{bc}] = -{\mathbb J}_{ac}\quad ,  \quad
[{\mathbb J}_{ac},{\mathbb J}_{bc}] = \omega_{bc}{\mathbb J}_{ab}\quad,
\label{comaa}
\end{equation}
where $\omega_{ab}=\prod_{s=a+1}^b\omega_s$ and $a<b<c$.
By simple rescaling of the generators the values $\omega_i$ may be brought to
one of the values 1, 0 or --1.

The structure of these algebras may be defined by two main statements:

$\bullet$
When all $\omega_i$ are non-zero the algebra is isomorphic to a certain
(pseudo-)orthogonal algebra.

$\bullet$
When a constant $\omega_a=0$ the resulting algebra
$so_{\omega_1,\dots,\omega_a=0,\dots,\omega_N}(N+1)$
has the semidirect structure
\begin{equation}
so_{\omega_1,\dots,\omega_a=0,\dots,\omega_N}(N+1)
\equiv t \odot
(so_{\omega_1,\dots,\omega_{a-1}}(a) \oplus
so_{\omega_{a+1},\dots,\omega_N}(N+1-a))
\quad,
\label{semidirect}
\end{equation}
where $t$ is an abelian subalgebra of dimension
$\text{dim}\,t=a(N+1-a)$ and  the remaining subalgebra is a direct sum.
In particular, when $a=1$ we obtain the usual (pseudo) orthogonal
inhomogeneous algebras
$so_{\omega_1=0,\omega_2,\dots,\omega_N}(N+1)$
with semidirect structure
\begin{equation}
so_{\omega_1=0,\omega_2,\dots,\omega_N}(N+1) \equiv
iso_{\omega_2,\dots,\omega_N}(N) =
t_N \odot so_{\omega_2,\dots,\omega_N}(N) \quad.
\label{semidirectin}
\end{equation}

The structure behind the decomposition (\ref{semidirect}) can be described 
visually by setting the generators in a triangular array (see Fig.~\ref{figure1}).
The generators spanning the subspace $t$ are those inside the rectangle,
while the subalgebras $so_{\omega_1,\dots,\omega_{a-1}}(a)$ and
$so_{\omega_{a+1},\dots,\omega_N}(N+1-a)$ correspond to the two triangles to
the left and below the rectangle respectively.
\begin{figure}[t]
\begin{center}
\begin{tabular}{cccc|cccc}
${\mathbb J}_{01} $&$ {\mathbb J}_{02} $&$\ldots$&
${\mathbb J}_{0\, a-1} $&
 ${\mathbb J}_{0a}$&${\mathbb J}_{0\, a+1}$&
$\ldots$&${\mathbb J}_{0N}$\\
 &$ {\mathbb J}_{12} $&$\ldots$&${\mathbb J}_{1\, a-1} $& ${\mathbb
J}_{1a}$&${\mathbb J}_{1 a+1}$&
$\ldots$&${\mathbb J}_{1N}$\\
 &&$\ddots $&$\vdots$&  $\vdots$&$\vdots$&
$ $&$\vdots$\\
 & &$ $&${\mathbb J}_{a-2\,a-1}$&
${\mathbb J}_{a-2\,a}$&${\mathbb J}_{a-2\,a+1}$&
$\ldots$&${\mathbb J}_{a-2\,N}$\\
 & &$ $& &  ${\mathbb J}_{a-1\,a}$&${\mathbb J}_{a-1\,a+1}$&
$\ldots$&${\mathbb J}_{a-1\,N}$\\
\cline{5-8}
 & &$ $&\multicolumn{1}{c}{\,}&    $ $&${\mathbb J}_{a\,a+1}$&
$\ldots$&${\mathbb J}_{a  N}$\\
 & &$ $&\multicolumn{1}{c}{\,}& $ $ &$ $&
$\ddots$&$\vdots$\\
 & &$ $&\multicolumn{1}{c}{\,}& $ $&$ $&
$ $&${\mathbb J}_{N-1\,N}$\\
\end{tabular}
\caption{Generators of the CK $so_{\omega_1,\dots,\omega_N}(N+1)$ algebra}
\label{figure1}
\end{center}
\end{figure}
In the $\omega_1=0$ ($\omega_N=0$) case the box is reduced to a single row
(column) in the large triangle.

To distinguish between the generators we shall denote by ${\mathbb X}$ those
inside the box (abelian algebra) and by ${\mathbb J}$ those in the two
triangles.
Namely,
\begin{equation}
\begin{array}{l}
{\mathbb X}_{ij} \Rightarrow i<a\ \text{and}\ j\ge a \quad,
\\
{\mathbb J}_{ij} \Rightarrow i\ge a\ \text{or}\ j< a \quad.
\end{array}
\label{defin}
\end{equation}
When two constants are set equal to zero ($\omega_a$ and $\omega_b$ say)
we have two different semidirect decompositions (\ref{semidirect})
corresponding to the constant $\omega_a=0$ or to the constant $\omega_b=0$.
For instance, the (3,1)--Galilei algebra appears in this context for
$\omega_1=0,\ \omega_2=0,\ \omega_3=1,\ \omega_4=1$ and accordingly has two
different semidirect structures which correspond to the constants $\omega_1$
and $\omega_2$.
In the triangular array this may be seen in Fig.~\ref{figure2} (for a
discussion on the dimensional analysis of the different contractions see
\cite[Section 2.2]{Azc.Olm.Bue.San:97}).

\begin{figure}[t]
\begin{center}
\begin{tabular}{|c|ccc}
\rule[-3pt]{0pt}{18pt}
$H$ & $P_1$ & $P_2$ & $P_3$
\\[0.3cm]
\cline{1-4}
\multicolumn{1}{c|}\
\rule[-3pt]{0pt}{18pt}&
$V_1$ & $V_2$ & $V_3$
\\[0.3cm]
\cline{2-4}
\multicolumn{1}{c}\
\rule[-3pt]{0pt}{18pt}
& & ${J}_3$ & $-{J}_2$ \\[0.3cm]
\multicolumn{1}{c}\
& & &\rule[-3pt]{0pt}{18pt}${J}_1$ \\[0.3cm]
\end{tabular}
\quad,\quad
\begin{tabular}{|c|ccc}
$T^{-1}$ & $T^{-1}TL^{-1}$ & $T^{-1}TL^{-1}$ & $T^{-1}TL^{-1}$
\\[0.3cm]
\cline{1-4}
\multicolumn{1}{c|}\
&
\rule[-3pt]{0pt}{18pt}
$TL^{-1}$ & $TL^{-1}$ & $TL^{-1}$
\\[0.3cm]
\cline{2-4}
\multicolumn{1}{c}\ \rule[-3pt]{0pt}{18pt}
& & 1 & 1 \\[0.3cm]
\multicolumn{1}{c}\ \rule[-3pt]{0pt}{18pt}
& & & 1 \\[0.3cm]
\end{tabular}
\caption{Generators of the Galilei algebra and its dimensional assignment}
\label{figure2}
\end{center}
\end{figure}

\section{The deformed family of inhomogeneous CK algebras}
\label{sec3}

Let us start with the set of inhomogeneous CK algebras $iso_{\omega_2\dots
\omega_N}(N)$ (see (\ref{semidirectin})).
There exists \cite{Bal.Her.Olm.San:94,Bal.Her.Olm.San:95}
a family of Hopf algebras, denoted by
${\cal U}_\lambda(iso_{\omega_2\dots \omega_N}(N))$, that are a deformation
of these CK algebras and, therefore, may be called `quantum' inhomogeneous CK
algebras.
In \cite{Azc.Olm.Bue.San:97} it was shown that all these deformed algebras are
endowed with a bicrossproduct structure that corresponds to the
undeformed semidirect one (\ref{semidirectin}) in which the abelian algebra is
given by the single row with generators ${\mathbb J}_{0i}$ (see
Fig.~\ref{figure1} for $a=1$).

Explicitly the deformed Hopf algebra
${\cal U}_\lambda(iso_{\omega_2\dots \omega_N}(N))$ is given (in the basis in
which its bicrossproduct structure is displayed) by

$\bullet$ Commutators
\begin{equation}
\begin{array}{l}
[{\mathbb J}_{0i},{\mathbb J}_{0j}]=0 \quad,
\quad [{\mathbb J}_{0i},{\mathbb J}_{0N}]=0\quad,
\\{}
[{\mathbb J}_{ij},{\mathbb J}_{ik}] =  \omega_{ij} {\mathbb J}_{jk}
\quad,\quad
   [{\mathbb J}_{ij},{\mathbb J}_{jk}] = -{\mathbb J}_{ik}   \quad,\quad
   [{\mathbb J}_{ik},{\mathbb J}_{jk}] = \omega_{jk}{\mathbb J}_{ij} \quad,
\\{}
[{\mathbb J}_{ij},{\mathbb J}_{iN}] =  \omega_{ij} {\mathbb J}_{jN}\quad,
\quad
   [{\mathbb J}_{ij},{\mathbb J}_{jN}] = -{\mathbb J}_{iN} \quad,  \quad
   [{\mathbb J}_{iN},{\mathbb J}_{jN}] = \omega_{jN}{\mathbb J}_{ij}\quad,
\\{}
[{\mathbb J}_{ij},{\mathbb J}_{0k}]=
\delta_{ik}{\mathbb J}_{0j} - \delta_{jk}\omega_{ij}{\mathbb
J}_{0i}\quad,\quad
   [{\mathbb J}_{ij},{\mathbb J}_{0N}]=0\quad,
\\{}
\displaystyle
[{\mathbb J}_{iN},{\mathbb J}_{0j}]=\delta_{ij}
\left(\frac{1-e^{-2\lambda {\mathbb J}_{0N}}}{2\lambda}-
     \frac \lambda 2  \sum_{s=1}^{N-1}\omega_{sN}{\mathbb J}_{0s}^2 \right)
     + \lambda \omega_{iN}{\mathbb J}_{0i}{\mathbb J}_{0j}\quad,
\\{}
[{\mathbb J}_{iN},{\mathbb J}_{0N}]=-\omega_{iN}{\mathbb J}_{0i} \quad .
\end{array}
\label{ce}
\end{equation}

$\bullet$ Coproduct
\begin{equation}
\begin{array}{l}
\Delta({\mathbb J}_{0i})=e^{-\lambda {\mathbb J}_{0N}}\otimes {\mathbb J}_{0i}
+ {\mathbb J}_{0i}\otimes  1 \quad, \quad
\Delta({{\mathbb J}_{0N}})=1\otimes{{\mathbb J}_{0N}}+{{\mathbb
J}_{0N}}\otimes
1 \quad,
\\[0.3cm]
\Delta({ {\mathbb J}_{ij}})=1\otimes{ {\mathbb J}_{ij}}+{{\mathbb
J}_{ij}}\otimes 1\quad,\\[0.3cm]
\displaystyle
\Delta({\mathbb J}_{iN})=e^{-\lambda {\mathbb J}_{0N}}\otimes {\mathbb J}_{iN}
+{\mathbb J}_{iN}\otimes 1
  +\lambda\sum_{s=1}^{i-1}\omega_{iN} {\mathbb J}_{0s} \otimes {\mathbb
J}_{si}
  -\lambda\sum_{s=i+1}^{N-1}\omega_{sN} {\mathbb J}_{0s} \otimes {\mathbb
J}_{is}\quad.
\label{cb}
\end{array}
\end{equation}

$\bullet$ Counit
\begin{equation}
\varepsilon({\mathbb J}_{0i})=\varepsilon({\mathbb J}_{0N})=
\varepsilon({\mathbb J}_{ij})=
\varepsilon({\mathbb J}_{iN})=0\quad.
\label{cc}
\end{equation}

$\bullet$ Antipode
\begin{equation}
\begin{array}{l}
\gamma({\mathbb J}_{0i})=-e^{\lambda {\mathbb J}_{0N}}
{\mathbb J}_{0i}\quad,\quad
\gamma({\mathbb J}_{0N})=-{\mathbb J}_{0N} \quad,\quad
\gamma({\mathbb J}_{ij})=-{\mathbb J}_{ij}\quad,
  \\[0.3cm]
\displaystyle
\gamma({\mathbb J}_{iN}) = - e^{\lambda {\mathbb J}_{0N}}{\mathbb J}_{iN} +
 \lambda e^{\lambda {\mathbb J}_{0N}}\sum_{s=1}^{i-1}\omega_{iN}
{\mathbb J}_{0s} {\mathbb J}_{si} -
 \lambda e^{\lambda {\mathbb J}_{0N}}\sum_{s=i+1}^{N-1}\omega_{sN}
{\mathbb J}_{0s} {\mathbb J}_{is} \quad.
\label{cd}
\end{array}
\end{equation}
In this basis it is easy to check the following
\begin{theorem}[{\rm \cite{Azc.Olm.Bue.San:97}}{}]
\label{th3.1}
The deformed Hopf CK family of algebras
${\cal U}_\lambda(iso_{\omega_2\dots \omega_N}(N))$
has a bicrossproduct structure
\begin{equation}
{\cal U}_\lambda(iso_{\omega_2\dots \omega_N}(N)) =
{\cal U}({so}_{\omega_2\dots \omega_N} (N))
^\beta\triangleright\!\!\!\blacktriangleleft_\alpha
{\cal U}_\lambda(T_N)
\label{biccd}
\end{equation}
relative to the right action
\begin{equation}
\alpha({\mathbb J}_{0i},{\mathbb J}_{jk})\equiv
{\mathbb J}_{0i} \triangleleft {\mathbb J}_{jk} :=
[{\mathbb J}_{0i},{\mathbb J}_{jk}]
\end{equation}
and left coaction $\beta$
\begin{equation}
\begin{array}{l}
\beta({\mathbb J}_{ij}) = 1\otimes {\mathbb J}_{ij} \quad,
\\[0.3cm]
\displaystyle
\beta({\mathbb J}_{iN}) :=e^{-\lambda {\mathbb J}_{0N}}\otimes
{\mathbb J}_{iN} +
\lambda\sum_{s=1}^{i-1}\omega_{iN} {\mathbb J}_{0s} \otimes {\mathbb J}_{si}
-\lambda\sum_{s=i+1}^{N-1}\omega_{sN} {\mathbb J}_{0s} \otimes {\mathbb
J}_{is}
\quad,
\end{array}
\end{equation}
where ${\cal U}_\lambda(T_N)$ is the abelian Hopf algebra generated by
${\mathbb J}_{0i}$ and ${\cal U}({so}_{\omega_2\dots \omega_N} (N))$ is the
\emph{undeformed} cocommutative Hopf algebra (with primitive coproduct)
generated by
${\mathbb J}_{ij}$ with the commutation relations given in the second and
third line of (\ref{ce}).
\end{theorem}

Let us now set $\omega_a=0$; then the algebra
${\cal U}_\lambda(iso_{\omega_2\dots \omega_a=0\dots \omega_N}(N))$ is given
(with the notation in (\ref{defin})) by

$\bullet$ Commutators
\begin{equation}
{\mathbb X}-\text{sector}
\left\{
[{\mathbb X}_{ij},{\mathbb X}_{kl}]=0
\rule[-6pt]{0pt}{18pt}
\right.
\label{Xsector}
\end{equation}
\begin{equation}
{\mathbb J}-\text{sector}
\left\{
\begin{array}{l}
[{\mathbb J}_{0i},{\mathbb J}_{0k}]=0
\\{}
[{\mathbb J}_{ij},{\mathbb J}_{0k}]=\delta_{ik}{\mathbb J}_{0j} -
\delta_{jk}\omega_{ij} {\mathbb J}_{0i}
\\{}
[{\mathbb J}_{iN},{\mathbb J}_{0j}]=\lambda \omega_{iN}
{\mathbb X}_{0i}{\mathbb J}_{0j}
\\{}
[{\mathbb J}_{ij},{\mathbb J}_{ik}]=\omega_{ij} {\mathbb J}_{jk}\ ,\
[{\mathbb J}_{ij},{\mathbb J}_{jk}]=-{\mathbb J}_{ik}\ ,\
[{\mathbb J}_{ik},{\mathbb J}_{jk}]=\omega_{jk} {\mathbb J}_{ij}
\\{}
[{\mathbb J}_{ij},{\mathbb J}_{iN}]=\omega_{ij} {\mathbb J}_{jN}\ ,\
[{\mathbb J}_{ij},{\mathbb J}_{jN}]=-{\mathbb J}_{iN}\ ,\
[{\mathbb J}_{iN},{\mathbb J}_{jN}]=\omega_{jN} {\mathbb J}_{ij}
\end{array}
\right.
\label{Jsector}
\end{equation}
\begin{equation}
{\mathbb J}{\mathbb X}-\text{sector}
\left\{
\begin{array}{l}
[{\mathbb J}_{0i},{\mathbb X}_{0j}]=[{\mathbb J}_{0i},{\mathbb X}_{0N}]=0
\\{}
[{\mathbb J}_{ij},{\mathbb X}_{ik}]=\omega_{ij} {\mathbb X}_{jk}\ ,\
[{\mathbb J}_{ij},{\mathbb X}_{jk}]=-{\mathbb X}_{ik}
\\{}
[{\mathbb J}_{jk},{\mathbb X}_{ij}]={\mathbb X}_{ik}\ ,\
[{\mathbb J}_{jk},{\mathbb X}_{ik}]=-\omega_{jk} {\mathbb X}_{ij}
\\{}
[{\mathbb J}_{ij},{\mathbb X}_{iN}]=\omega_{ij} {\mathbb X}_{jN}\ ,\
[{\mathbb J}_{ij},{\mathbb X}_{jN}]=-{\mathbb X}_{iN}
\\{}
[{\mathbb J}_{jN},{\mathbb X}_{ij}]={\mathbb X}_{iN}\ ,\
[{\mathbb J}_{jN},{\mathbb X}_{iN}]=-\omega_{jN} {\mathbb X}_{ij}
\\{}
[{\mathbb J}_{0k},{\mathbb X}_{ij}]=-\delta_{ik}{\mathbb X}_{0j}
\ ,\
[{\mathbb J}_{ij},{\mathbb X}_{0k}]=
\delta_{ik}{\mathbb X}_{0j} - \delta_{jk}\omega_{ij} {\mathbb X}_{0i}
\\{}
[{\mathbb J}_{ij},{\mathbb X}_{0N}]=0
\ ,\
[{\mathbb J}_{iN},{\mathbb X}_{0N}]=-\omega_{iN} {\mathbb X}_{0i}
\\{}
\displaystyle
[{\mathbb J}_{0i},{\mathbb X}_{jN}]=
-\delta_{ij}\left( {1-e^{-2\lambda{\mathbb X}_{0N}}\over2\lambda}
-{\lambda\over 2}\sum_{s=a}^{N-1}\omega_{sN} {\mathbb X}_{0s}^2 \right)
\\{}
\displaystyle
[{\mathbb J}_{iN},{\mathbb X}_{0j}]=
\delta_{ij}\left( {1-e^{-2\lambda{\mathbb X}_{0N}}\over2\lambda}
-{\lambda\over 2}\sum_{s=a}^{N-1}\omega_{sN} {\mathbb X}_{0s}^2 \right)
+\lambda\omega_{iN}{\mathbb X}_{0i}{\mathbb X}_{0j}
\end{array}
\right.
\label{JXsector}
\end{equation}

$\bullet$ Coproduct
\begin{equation}
{\mathbb X}-\text{sector}
\left\{
\begin{array}{l}
\Delta {\mathbb X}_{0N} = 1\otimes {\mathbb X}_{0N} + {\mathbb X}_{0N}\otimes
1
\ ,\
\Delta {\mathbb X}_{ij} = 1\otimes {\mathbb X}_{ij} + {\mathbb X}_{ij}\otimes
1
\\{}
\Delta {\mathbb X}_{0i} = e^{-\lambda {\mathbb X}_{0N}}
\otimes {\mathbb X}_{0i} + {\mathbb X}_{0i}\otimes 1
\\{}
\displaystyle
\Delta {\mathbb X}_{iN}=
e^{-\lambda {\mathbb X}_{0N}} \otimes {\mathbb X}_{iN} +
{\mathbb X}_{iN}\otimes 1
-\lambda \sum_{s=a}^{N-1}\omega_{sN} {\mathbb X}_{0s} \otimes {\mathbb X}_{is}
\end{array}
\right.
\label{coXsec}
\end{equation}

\begin{equation}
{\mathbb J}-\text{sector}
\left\{
\begin{array}{l}
\Delta {\mathbb J}_{0i} = e^{-\lambda {\mathbb X}_{0N}}
\otimes {\mathbb J}_{0i} + {\mathbb J}_{0i}\otimes 1 \ ,\
\Delta {\mathbb J}_{ij} = 1\otimes {\mathbb J}_{ij} + {\mathbb J}_{ij}\otimes
1
\\{}
\displaystyle
\Delta {\mathbb J}_{iN}=
e^{-\lambda {\mathbb X}_{0N}} \otimes {\mathbb J}_{iN}
+{\mathbb J}_{iN}\otimes 1
+\lambda \sum_{s=1}^{a-1}\omega_{iN} {\mathbb J}_{0s} \otimes {\mathbb X}_{si}
\\{}
\displaystyle
\qquad\qquad
+\lambda \sum_{s=a}^{i-1}\omega_{iN} {\mathbb X}_{0s} \otimes {\mathbb J}_{si}
-\lambda \sum_{s=i+1}^{N-1}\omega_{sN} {\mathbb X}_{0s} \otimes {\mathbb
J}_{is}
\end{array}
\right.
\label{coJsec}
\end{equation}

$\bullet$ Counit
\begin{equation}
\begin{array}{l}
\varepsilon({\mathbb J}_{ij})=\varepsilon({\mathbb J}_{0i})=
\varepsilon({\mathbb J}_{iN})=0
\\
\varepsilon({\mathbb X}_{ij})=\varepsilon({\mathbb X}_{0i})=
\varepsilon({\mathbb X}_{0N})=\varepsilon({\mathbb X}_{iN})=0
\end{array}
\label{coun}
\end{equation}

$\bullet$ Antipode
\begin{equation}
{\mathbb X}-\text{sector}
\left\{
\begin{array}{l}
\gamma({\mathbb X}_{0N})=-{\mathbb X}_{0N}
\ ,\
\gamma({\mathbb X}_{ij})=-{\mathbb X}_{ij}
\\
\gamma({\mathbb X}_{0i})=-e^{\lambda{\mathbb X}_{0N}}{\mathbb X}_{0i}
\\
\displaystyle
\gamma({\mathbb X}_{iN})=-e^{\lambda{\mathbb X}_{0N}}{\mathbb X}_{iN}
-\lambda e^{\lambda{\mathbb X}_{0N}}
\sum_{s=a}^{N-1}\omega_{sN}{\mathbb X}_{0s}{\mathbb X}_{is}
\end{array}
\right.
\label{Xantipod}
\end{equation}

\begin{equation}
{\mathbb J}-\text{sector}
\left\{
\begin{array}{l}
\gamma({\mathbb J}_{ij})=-{\mathbb J}_{ij}
\\
\gamma({\mathbb J}_{0i})=-e^{\lambda{\mathbb X}_{0N}}{\mathbb J}_{0i}
\\
\displaystyle
\gamma({\mathbb J}_{iN})=-e^{\lambda{\mathbb X}_{0N}}{\mathbb J}_{iN}
+\lambda e^{\lambda{\mathbb X}_{0N}}
\sum_{s=1}^{a-1}\omega_{iN} {\mathbb J}_{0s}{\mathbb X}_{si}
\\
\qquad\qquad
\displaystyle
+\lambda e^{\lambda{\mathbb X}_{0N}}
\sum_{s=a}^{i-1}\omega_{iN} {\mathbb X}_{0s}{\mathbb J}_{si}
-\lambda e^{\lambda{\mathbb X}_{0N}}\sum_{s=i+1}^{N-1}{\mathbb X}_{0s}{\mathbb
J}_{is}
\end{array}
\right.
\label{Jantipod}
\end{equation}
This algebra, as result of theorem~\ref{th3.1}, has a bicrossproduct
structure (\ref{biccd}).
However the theorem does not give us information about the (possible)
bicrossproduct structure for the decomposition given in (\ref{semidirect}).
This is the problem that we address now.

\section{Bicrossproduct structure}
\label{sec4}

The algebra given above (\ref{Xsector})-(\ref{Jantipod})
does not present directly a bicrossproduct structure for the decomposition
(\ref{semidirect}).
This is due to the term
$\lambda \sum_{s=1}^{a-1}\omega_{iN} {\mathbb J}_{0s} \otimes {\mathbb
X}_{si}$
in the ${\mathbb J}_{iN}$ coproduct (second line in (\ref{coJsec}))
and to the commutator
$[{\mathbb J}_{iN},{\mathbb J}_{0j}]$ (third line in (\ref{Jsector})) that
does not close a ${\mathbb J}$ algebra\footnote{\label{foot4a}
Nevertheless in the particular
case $a=N$ ($\omega_N=0$) these terms are not present, and the change of basis
given in (\ref{newset}) is not necessary (notice that for $\omega_N=0$ the
change of basis is trivial).}.

Let us define
\begin{equation}
\hat J_{iN} =
\lambda \sum_{s=1}^{a-1}\omega_{iN} {\mathbb J}_{0s} {\mathbb X}_{si}
\label{hatJ}
\end{equation}
that verifies
\begin{equation}
\Delta \hat J_{iN} = e^{-\lambda {\mathbb X}_{0N}}\otimes \hat J_{iN}
+  \hat J_{iN}\otimes 1 +
\lambda \sum_{s=1}^{a-1}\omega_{iN} {\mathbb J}_{0s} \otimes {\mathbb X}_{si}
+ \lambda \sum_{s=1}^{a-1}\omega_{iN} e^{-\lambda {\mathbb X}_{0N}}
{\mathbb X}_{si}\otimes {\mathbb J}_{0s}
\label{dem1}
\end{equation}
and
\begin{equation}
[\hat J_{iN}, {\mathbb J}_{0j}]=\lambda\omega_{iN} {\mathbb J}_{0j}
{\mathbb X}_{0i}\quad.
\label{dem2}
\end{equation}
Thus, the change of basis
${\mathbb J}_{iN} \to {\mathbb J}_{iN} - \hat J_{iN}$
solves the two difficulties pointed out before.

Specifically, if we introduce the new set of generators
\begin{equation}
\begin{array}{c}
J_{0i}= {\mathbb J}_{0i}\quad,\quad
X_{0i}= {\mathbb X}_{0i}\quad,\quad
X_{0N}= {\mathbb X}_{0N}\quad,
\\
J_{ij}= {\mathbb J}_{ij}\quad,\quad
X_{ij}= {\mathbb X}_{ij}\quad,\quad
X_{iN}= {\mathbb X}_{iN} \quad,
\\
J_{iN} = {\mathbb J}_{iN} - \hat J_{iN}
\end{array}
\label{newset}
\end{equation}
the algebra
${\cal U}_\lambda(iso_{\omega_2,\dots,\omega_a=0,\dots, \omega_N}(N))$ is
written as

$\bullet$ Commutators
\begin{equation}
{X}-\text{sector}
\left\{
\begin{array}{l}
[X_{0i},X_{0j}]=[X_{0i},X_{0N}]=0
\\{}
[X_{ij},X_{0k}]=[{X}_{ij},{X}_{kl}]=[X_{ij},X_{kN}]=0
\\{}
[X_{iN},X_{0j}]=[X_{iN},X_{0N}]=[X_{iN},X_{jk}]=[X_{iN},X_{jN}]=0
\end{array}
\right.
\label{X-sector}
\end{equation}
\begin{equation}
{J}-\text{sector}
\left\{
\begin{array}{l}
[{J}_{0i},{J}_{0k}]=0
\\{}
[{J}_{ij},{J}_{0k}]=\delta_{ik}{J}_{0j} -
\delta_{jk}\omega_{ij} J_{0i}
\\{}
[{J}_{iN},{J}_{0j}]=0
\\{}
[{J}_{ij},{J}_{ik}]=\omega_{ij} {J}_{jk}\ ,\
[{J}_{ij},{J}_{jk}]=-{J}_{ik}\ ,\
[{J}_{ik},{J}_{jk}]=\omega_{jk} {J}_{ij}
\\{}
[{J}_{ij},{J}_{iN}]=\omega_{ij} {J}_{jN}\ ,\
[{J}_{ij},{J}_{jN}]=-{J}_{iN}\ ,\
[{J}_{iN},{J}_{jN}]=\omega_{jN} {J}_{ij}
\end{array}
\right.
\label{J-sector}
\end{equation}
\begin{equation}
{JX}-\text{sector}
\left\{
\begin{array}{l}
[{J}_{0i},{X}_{0j}]=[{J}_{0i},{X}_{0N}]=0
\\{}
[{J}_{ij},{X}_{ik}]=\omega_{ij} {X}_{jk}\ ,\
[{J}_{ij},{X}_{jk}]=-{X}_{ik}
\\{}
[{J}_{jk},{X}_{ij}]={X}_{ik}\ ,\
[{J}_{jk},{X}_{ik}]=-\omega_{jk} {X}_{ij}
\\{}
[{J}_{ij},{X}_{iN}]=\omega_{ij} {X}_{jN}\ ,\
[{J}_{ij},{X}_{jN}]=-{X}_{iN}\ ,\
[{J}_{jN},{X}_{ij}]={X}_{iN}
\\{}
\displaystyle
[{J}_{jN},{X}_{iN}]=-\omega_{jN} {X}_{ij}
+\lambda \omega_{jN} \left( {1-e^{-2\lambda{X}_{0N}}\over2\lambda}
-{\lambda\over 2}\sum_{s=a}^{N-1}\omega_{sN} {X}_{0s}^2 \right) X_{ij}
\\{}
[{J}_{0k},{X}_{ij}]=-\delta_{ik}{X}_{0j}
\ ,\
[{J}_{ij},{X}_{0k}]=
\delta_{ik}{X}_{0j} - \delta_{jk}\omega_{ij} {X}_{0i}
\\{}
[{J}_{ij},{X}_{0N}]=0
\ ,\
[{J}_{iN},{X}_{0N}]=-\omega_{iN} {X}_{0i}
\\{}
[J_{kN},X_{ij}]=\delta_{jk} X_{iN} + \lambda \omega_{kN} X_{0j} X_{ik}
\\{}
\displaystyle
[{J}_{0i},{X}_{jN}]=
-\delta_{ij}\left( {1-e^{-2\lambda{X}_{0N}}\over2\lambda}
-{\lambda\over 2}\sum_{s=a}^{N-1}\omega_{sN} {X}_{0s}^2 \right)
\\{}
\displaystyle
[{J}_{iN},{X}_{0j}]=
\delta_{ij}\left( {1-e^{-2\lambda{X}_{0N}}\over2\lambda}
-{\lambda\over 2}\sum_{s=a}^{N-1}\omega_{sN} {X}_{0s}^2 \right)
+\lambda\omega_{iN}{X}_{0i}{X}_{0j}
\end{array}
\right.
\label{JX-sector}
\end{equation}

$\bullet$ Coproduct
\begin{equation}
{X}-\text{sector}
\left\{
\begin{array}{l}
\Delta {X}_{0N} = 1\otimes {X}_{0N} + {X}_{0N}\otimes 1
\ ,\
\Delta {X}_{ij} = 1\otimes {X}_{ij} + {X}_{ij}\otimes 1
\\{}
\Delta {X}_{0i} = e^{-\lambda {X}_{0N}}
\otimes {X}_{0i} + {X}_{0i}\otimes 1
\\{}
\displaystyle
\Delta {X}_{iN}=
e^{-\lambda {X}_{0N}} \otimes {X}_{iN}+
{X}_{iN}\otimes 1
-\lambda \sum_{s=a}^{N-1}\omega_{sN} {X}_{0s} \otimes {X}_{is}
\end{array}
\right.
\label{coX-sec}
\end{equation}

\begin{equation}
{J}-\text{sector}
\left\{
\begin{array}{l}
\Delta {J}_{0i} = e^{-\lambda {X}_{0N}}
\otimes {J}_{0i} + {J}_{0i}\otimes 1 \ ,\
\Delta {J}_{ij} = 1\otimes {J}_{ij} + {J}_{ij}\otimes 1
\\{}
\displaystyle
\Delta {J}_{iN}=
e^{-\lambda {X}_{0N}} \otimes {J}_{iN}
+{J}_{iN}\otimes 1
-\lambda \sum_{s=1}^{a-1}\omega_{iN} e^{-\lambda {X}_{0N}}
{X}_{si} \otimes {J}_{0s}
\\{}
\displaystyle
\qquad\qquad
+\lambda \sum_{s=a}^{i-1}\omega_{iN} {X}_{0s} \otimes {J}_{si}
-\lambda \sum_{s=i+1}^{N-1}\omega_{sN} {X}_{0s} \otimes {J}_{is}
\end{array}
\right.
\label{coJ-sec}
\end{equation}

$\bullet$ Counit
\begin{equation}
\begin{array}{l}
\varepsilon({J}_{ij})=\varepsilon({J}_{0i})=
\varepsilon({J}_{iN})=0
\\
\varepsilon({X}_{ij})=\varepsilon({X}_{0i})=\varepsilon({X}_{0N})=
\varepsilon({X}_{iN})=0
\end{array}
\label{coun-sec}
\end{equation}

$\bullet$ Antipode
\begin{equation}
{X}-\text{sector}
\left\{
\begin{array}{l}
\gamma({X}_{0N})=-{X}_{0N}
\ ,\
\gamma({X}_{ij})=-{X}_{ij}
\\
\gamma({X}_{0i})=-e^{\lambda{X}_{0N}}{X}_{0i}
\\
\displaystyle
\gamma({X}_{iN})=-e^{\lambda{X}_{0N}}{X}_{iN}
-\lambda e^{\lambda{X}_{0N}}\sum_{s=a}^{N-1}\omega_{sN}{X}_{0s}{X}_{is}
\end{array}
\right.
\label{X-antipod}
\end{equation}

\begin{equation}
{J}-\text{sector}
\left\{
\begin{array}{l}
\gamma({J}_{ij})=-{J}_{ij}
\\
\gamma({J}_{0i})=-e^
{\lambda{X}_{0N}}{J}_{0i}
\\
\displaystyle
\gamma({J}_{iN})=-e^{\lambda{X}_{0N}}{J}_{iN}
-\lambda e^{\lambda{X}_{0N}}\sum_{s=1}^{a-1}\omega_{iN} {X}_{si}{J}_{0s}
\\
\qquad\qquad
\displaystyle
+\lambda e^{\lambda{X}_{0N}}\sum_{s=a}^{i-1}\omega_{iN} {X}_{0s}{J}_{si}
-\lambda e^{\lambda{X}_{0N}}\sum_{s=i+1}^{N-1} \omega_{sN} {X}_{0s}{X}_{is}
\end{array}
\right.
\label{J-antipod}
\end{equation}

In this basis we now state the following
\begin{theorem}
\label{theor1}
The algebra
${\cal U}_\lambda(iso_{\omega_2\dots \omega_a=0\dots \omega_N}(N))$
has the bicrossproduct structure
\begin{equation}
{\cal U}_\lambda(iso_{\omega_2\dots \omega_a=0\dots \omega_N}(N))
=
{\cal U}(so_{\omega_1=0,\dots,\omega_{a-1}}(a))
\oplus
{\cal U}(so_{\omega_{a+1},\dots,\omega_N}(N+1-a))
^\beta\triangleright\!\!\!\blacktriangleleft_\alpha
{\cal U}_\lambda(T_{a(N+1-a)})
\label{theor1a}
\end{equation}
where ${\cal U}(so_{\omega_1=0,\dots,\omega_{a-1}}(a))$ is the undeformed Hopf
algebra generated by $\{J_{ij},\ i<(a-1),\ j<a\}$,
${\cal U}(so_{\omega_{a+1},\dots,\omega_N}(N+1-a))$ is spanned by the
generators $\{J_{ij},\ i>(a-1),\ j>a\}$
and ${\cal U}_\lambda(T_{a(N+1-a)})$ is the deformed abelian algebra generated
by $X_{ij}$ (recall that those generators are
restricted to the indices
$i<a\ \text{and}\ j\ge a$, (\ref{defin})).
The right action
$\alpha:
{\cal U}_\lambda(T_{a(N+1-a)})
\otimes {\cal U}(so_{\omega_1=0,\dots,\omega_{a-1}}(a))
\oplus {\cal U}(so_{\omega_{a+1},\dots,\omega_N}(N+1-a))
\to {\cal U}_\lambda(T_{a(N+1-a)})$
is defined by (\ref{JX-sector}) through
\begin{equation}
\alpha(X_{ij},J_{kl})\equiv X_{ij}\triangleleft J_{kl} :=
[X_{ij},J_{kl}]
\label{alphadef}
\end{equation}
and the (left) coaction $\beta:
{\cal U}(so_{\omega_1=0,\dots,\omega_{a-1}}(a))
\oplus {\cal U}(so_{\omega_{a+1},\dots,\omega_N}(N+1-a))
\to {\cal U}_\lambda(T_{a(N+1-a)}) \otimes
{\cal U}(so_{\omega_1=0,\dots,\omega_{a-1}}(a))
\oplus {\cal U}(so_{\omega_{a+1},\dots,\omega_N}(N+1-a))$
is designed to reproduce the coproduct (\ref{coJ-sec})
\begin{equation}
\begin{array}{l}
\beta(J_{ij})=1\otimes J_{ij}
\\
\beta(J_{0i})=e^{-\lambda X_{0N}}\otimes J_{ij}
\\
\displaystyle
\beta(J_{iN})=e^{-\lambda X_{0N}}\otimes J_{iN}
-\lambda\sum_{s=1}^{a-1}\omega_{iN}e^{-\lambda X_{0N}} X_{si}\otimes J_{0s}
+\lambda\sum_{s=a}^{i-1}\omega_{iN} X_{0s}\otimes J_{si}
\\
\displaystyle
\qquad\qquad
-\lambda\sum_{s=i+1}^{N-1}\omega_{sN} X_{0s}\otimes J_{is}
\end{array}
\label{betadef}
\end{equation}
\end{theorem}

{}From theorem~\ref{theor1} it is easy to check the following

\begin{corollary}
\label{cor4.1}
Associated to each graded contraction in the inhomogeneous CK family of
deformed algebras we have a different bicrossproduct structure related to the
corresponding Lie algebra semidirect structure that appears in the contraction
(see (\ref{semidirect})).
These bicrossproduct structures are preserved under further (graded)
contraction processes.
\end{corollary}
\section{Examples}

\subsection{$N=3$ case}

In the $N=3$ case we obtain (in the basis of (\ref{ce})-(\ref{cd})) the
following equations
\begin{equation}
\begin{array}{l}
[{\mathbb J}_{0i},{\mathbb J}_{0j}]=0\quad,\quad i,j=1,2,3 \quad,
\\{}
[{\mathbb J}_{12},{\mathbb J}_{13}]=\omega_2 {\mathbb J}_{23} \quad,\quad
[{\mathbb J}_{13},{\mathbb J}_{23}]=\omega_3 {\mathbb J}_{12} \quad,\quad
[{\mathbb J}_{12},{\mathbb J}_{23}]=-{\mathbb J}_{13}\quad,
\\{}
[{\mathbb J}_{12},{\mathbb J}_{03}]=0\quad,\quad
[{\mathbb J}_{12},{\mathbb J}_{01}]={\mathbb J}_{02}\quad,\quad
[{\mathbb J}_{12},{\mathbb J}_{02}]=-\omega_2 {\mathbb J}_{01}\quad,
\\{}
[{\mathbb J}_{13},{\mathbb J}_{03}]=-\omega_2 \omega_3 {\mathbb J}_{01}
\quad,\quad
[{\mathbb J}_{23},{\mathbb J}_{03}]=-\omega_2 {\mathbb J}_{02}\quad,
\\
\displaystyle
[{\mathbb J}_{13},{\mathbb J}_{01}]=
{1-e^{-2\lambda {\mathbb J}_{03}}\over 2\lambda} +
{\lambda\over 2}\omega_3 (\omega_2 {\mathbb J}_{01}^2- {\mathbb J}_{02}^2)
\quad,
\\
\displaystyle
[{\mathbb J}_{23},{\mathbb J}_{02}]=
{1-e^{-2\lambda {\mathbb J}_{03}}\over 2\lambda} -{\lambda\over 2}\omega_3
(\omega_2 {\mathbb J}_{01}^2- {\mathbb J}_{02}^2)
\quad,
\\{}
[{\mathbb J}_{13},{\mathbb J}_{02}]=
\lambda \omega_2 \omega_3 {\mathbb J}_{01} {\mathbb J}_{02}
\quad,\quad
[{\mathbb J}_{23},{\mathbb J}_{01}]=
\fbox{ $\lambda \omega_3  {\mathbb J}_{01} {\mathbb J}_{02}$ } \quad,
\end{array}
\end{equation}
\begin{equation}
\begin{array}{l}
\Delta {\mathbb J}_{0i}=
e^{-\lambda {\mathbb J}_{03}}\otimes {\mathbb J}_{0i} +
{\mathbb J}_{0i} \otimes 1 \quad,\quad i=1,2\quad,
\\
\Delta {\mathbb J}_{03} = 1\otimes {\mathbb J}_{03} + {\mathbb J}_{03}
\otimes 1 \quad,\quad
\Delta {\mathbb J}_{12} = 1\otimes {\mathbb J}_{12} +
{\mathbb J}_{12} \otimes 1\quad,
\\
\Delta {\mathbb J}_{13} = e^{-\lambda {\mathbb J}_{03}}\otimes
{\mathbb J}_{13}  + {\mathbb J}_{13} \otimes 1
-\lambda \omega_3 {\mathbb J}_{02}\otimes {\mathbb J}_{12} \quad,
\\
\Delta {\mathbb J}_{23} = e^{-\lambda {\mathbb J}_{03}}\otimes
{\mathbb J}_{23}  + {\mathbb J}_{23} \otimes 1
+ \fbox{ $\lambda \omega_3 {\mathbb J}_{01}\otimes {\mathbb J}_{12}$ } \quad.
\end{array}
\end{equation}
We have stressed with a \fbox{\rule{10pt}{0pt}box\rule{10pt}{0pt}}
the terms that might not allow the algebra to be a
bicrossproduct (see first paragraph in sec. \ref{sec4}).
If $\omega_3=0$ these terms cancel\footnote{Note that for $\omega_3=0$ (\ie,
$\omega_N=0$) the change of basis (\ref{newset}) is trivial
(see (\ref{hatJ}) and footnote \ref{foot4a}).}
but if $\omega_2=0$ we keep them
and we need the change of basis given in (\ref{newset}).
In the new basis we have
\begin{equation}
\Delta J_{23} = e^{-\lambda X_{03}}\otimes J_{23}  + J_{23} \otimes 1
-\lambda \omega_3 X_{12} e^{-\lambda X_{03}}\otimes J_{01}
\end{equation}
and
\begin{equation}
[J_{23},J_{01}]=0 \quad.
\end{equation}
Thus, the terms marked above have disappeared after the change of basis and we
have a bicrossproduct structure as given in theorem~\ref{theor1}.

\subsection{A particular case: the Heisenberg-Weyl algebra}

Now we are going to study the case $a=N$
(\ie, $\omega_N=0$).
First, let us rename the generators in the basis (\ref{newset}) as
\begin{equation}
J_{0i}= X_{i} \quad,\quad
J_{ij}= J_{ij} \quad,\quad
X_{iN}= Y_{i} \quad,\quad
X_{0N}=\Xi
\label{recall}
\end{equation}
(note that for $\omega_N=0$ the $X$ sector is reduced to a single column in
the
triangular array in Fig.~\ref{figure1}).
Now the equations (\ref{X-sector}), (\ref{J-sector}) and (\ref{JX-sector})
acquire the form
\begin{equation}
\begin{array}{l}
[X_i,\Xi]=[Y_i,\Xi]=[J_{ij},\Xi]=0 \quad,
\\{}
[X_i,X_j]=0\quad,\quad [Y_i,Y_j]=0\quad,
\\{}
[J_{ij},J_{ik}]=\omega_{ij} J_{jk}\quad,\quad
[J_{ij},J_{jk}]=-J_{ik}\quad,\quad
[J_{ik},J_{jk}]=\omega_{jk} J_{ij}\quad,
\\{}
[J_{ij},X_{k}]=\delta_{ik}X_{j} -\delta_{jk}\omega_{ij}X_{i}
\quad,\quad
[J_{ij},Y_{k}]=\delta_{ik}\omega_{ij}Y_{j} -\delta_{jk}Y_{i} \quad,
\\{}
\displaystyle
[X_{i},Y_{j}]=-\delta_{ij}\left({1-e^{-2\lambda\Xi}\over 2\lambda} \right)
\quad.
\end{array}
\label{HWalgebra}
\end{equation}
In this way, we easily recognize the deformed Heisenberg-Weyl (HW) algebra
\cite{Azc.Bue:96} where $\Xi$ is the
central generator and the $J_{ij}$ generators act as a rotation group on the
$X_{i}$ and $Y_{i}$ generators.
The coproduct (\ref{coX-sec})-(\ref{coJ-sec}) takes the form
\begin{equation}
\begin{array}{l}
\Delta \Xi = 1\otimes \Xi + \Xi\otimes 1
\quad,\quad
\Delta J_{ij} =1\otimes J_{ij} + J_{ij}\otimes 1 \quad,
\\{}
\Delta X_i =e^{-\lambda\Xi} \otimes X_i + X_i\otimes 1
\quad,\quad
\Delta Y_i =e^{-\lambda\Xi} \otimes Y_i + Y_i\otimes 1 \quad.
\end{array}
\label{HWcoproduct}
\end{equation}
{}From the arguments given above we know that this algebra has two different
bicrossproduct (semidirect like) structures, one for the abelian algebra
generated by $\{X_i,\Xi\}$, and the other for the abelian algebra generated by
$\{Y_i,\Xi\}$ \cite{Azc.Bue:96}.

But, in this case, we have an additional cocycle-bicrossproduct structure 
(analogue to the undeformed central extension structure of the HW-algebra).
To see this let us define the algebra ${\cal H}$ as the undeformed algebra
generated by $\{X_i,Y_i,J_{ij}\}$ with primitive coproduct and commutators
\begin{equation}
\begin{array}{l}
[J_{ij},J_{ik}]=\omega_{ij} J_{jk}\quad,\quad
[J_{ij},J_{jk}]=-J_{ik}\quad,\quad
[J_{ik},J_{jk}]=\omega_{jk} J_{ij}\quad,
\\{}
[J_{ij},X_{k}]=\delta_{ik}X_{j} -\delta_{jk}\omega_{ij}X_{i}
\quad,\quad
[J_{ij},Y_{k}]=\delta_{ik}\omega_{ij}Y_{j} -\delta_{jk}Y_{i} \quad,
\\{}
[X_i,X_j]=0\quad,\quad [Y_i,Y_j]=0\quad,\quad [X_i,Y_j]=0
\end{array}
\label{HWalg1}
\end{equation}
(note that all the commutators are identical to those in (\ref{HWalgebra})
but the $[X_i,Y_j]$ one that now is abelian).
The algebra ${\cal A}$ is the undeformed algebra ${\cal U}(\Xi)$.
Now if we define the right action $\triangleleft : {\cal A}\otimes {\cal H}
\to {\cal A}$
\begin{equation}
\Xi \triangleleft J_{ij} = 0
\quad,\quad
\Xi \triangleleft X_i =0
\quad,\quad
\Xi \triangleleft Y_i =0
\label{HWright}
\end{equation}
(central extension means trivial action), the left coaction
$\beta:{\cal H}\to {\cal A}\otimes {\cal H}$
\begin{equation}
\beta (J_{ij})=1\otimes J_{ij}
\quad,\quad
\beta(X_i)=e^{-\lambda\Xi}\otimes X_i
\quad,\quad
\beta(Y_i)=e^{-\lambda\Xi}\otimes X_i\quad,
\label{HWleft}
\end{equation}
the antisymmetric two-cocycle $\xi:{\cal H}\otimes {\cal H} \to {\cal A}$
\footnote{\label{foot1}
The antisymmetric form of the cocycle is a matter of convention; different
forms of the cocycle are related by a coboundary change (see \cite{Azc.Bue:96b}
for
an explicit example).}
\begin{equation}
\xi (X_i, Y_j) = - \xi(Y_j, X_i)
= -{\delta_{ij}\over 2} \left({1-e^{-2\lambda\Xi}\over 2\lambda} \right)
\label{HWcocycle}
\end{equation}
and a trivial `two-cococycle' the HW algebra is given by the bicrossproduct
\begin{equation}
{\cal U}_\lambda (HW) = {\cal H}
\triangleright \!\!\! \blacktriangleleft_\xi {\cal A}\quad.
\end{equation}

In this form it is easy to recover the dual algebra
$\text{Fun}_\lambda(\text{HW})$ \cite{Azc.Bue:96}.
Let $R_{ij}$ be the dual generators corresponding to the undeformed `rotation'
algebra generated by $J_{ij}$~\footnote{This algebra is a true rotation
algebra for $\omega_i =1$ $i=1,\dots.N-1$; in general it is an inhomogeneous
algebra (if some $\omega=0$) or a pseudo-orthogonal algebra.
The dual algebra is given by the matrix representation $R_{ij}$.}
and let $x_i,\ y_j$ be the dual coordinates to the generators $X_i,\ Y_j$.
Then, the algebra $H$ dual to ${\cal H}$ is given by
\begin{equation}
\begin{array}{c}
\Delta R_{ij} = R_{ik} \otimes R_{kj}\quad,
\\{}
\Delta x_i = 1\otimes x_i + x_k\otimes R_{ki}
\quad,\quad
\Delta y_i = 1\otimes y_i + y_k \otimes R^{-1}_{ik}\quad;
\end{array}
\label{copHW}
\end{equation}
\begin{equation}
[R_{ij},R_{kl}]=[R_{ij},x_k]=[R_{ij},y_k]=[x_i,y_j]=0\quad.
\label{commHW}
\end{equation}
If we introduce the coordinate $\chi$ dual to the central generator $\Xi$ we
may
complete the dual algebra by dualizing the left coaction (\ref{HWleft}) and
the two-cocycle (\ref{HWcocycle}).
The left action is defined as the dual to the left coaction
\begin{equation}
\chi\triangleright x_i =[\chi,x_i]=-\lambda x_i
\quad,\quad
\chi\triangleright y_i =[\chi,y_i]=-\lambda y_i
\quad,\quad
\chi\triangleright R_{ij} =[\chi,R_{ij}]= 0
\label{rightHW}
\end{equation}
and the dual to the two-cocycle defines the two-cococycle\footnote{
As said in footnote \ref{foot1} we may choose a different form of the
two-cococycle.
For instance $\bar\psi(\chi)= y_i\otimes R_{ji}^{-1} x_j$ is also a
two-cococycle (related to (\ref{twococo}) by the cocoboundary
${1\over 2}y_i x_i$).}
\begin{equation}
\bar\psi (\chi) =
{1\over 2} ( y_i\otimes R_{ji}^{-1} x_j - x_i\otimes R_{ij} y_j )\quad.
\label{twococo}
\end{equation}
Thus, the coproduct is given by
\begin{equation}
\Delta \chi = 1\otimes\chi +\chi\otimes 1 +
{1\over 2} ( y_i\otimes R_{ji}^{-1} x_j - x_i\otimes R_{ij} y_j )\quad.
\end{equation}
As we may see the bicrossproduct structure (with cocycle in this case)
allows us to recover $\text{Fun}_\lambda(\text{HW})$ in an easy way from the
enveloping (dual) algebra ${\cal U}_\lambda(\text{HW})$.

\section*{Acknowledgements}
The author thanks J. A. de Azc\'arraga for his comments on the manuscript.
This paper has been partially supported by a research grant (PB096--0756)
from the MEC, Spain.
The author wishes to acknowledge an FPI grant from the Spanish Ministry of
Education and Culture and the CSIC.



\begin{thebibliography}{10}

\bibitem{Dri:87}
V.~G. Drinfel'd,
\newblock {\em Quantum groups},
\newblock in A.~Gleason, editor, {\em 1986 Int. Congr. of Math., MSRI}, vol.~I,
  page 798, American Mathematical Society, Providence, 1987.

\bibitem{Jim:85}
M.~Jimbo,
\newblock {\em A $q$-difference analogue of $U(g)$ and the Yang-Baxter
  equation},
\newblock Lett. Math. Phys. {\bf 10}, 63--69 (1985).

\bibitem{Jim:86}
M.~Jimbo,
\newblock {\em A $q$-analogue of $U(ql(N+1))$, Hecke algebra, and the
  Yang-Baxter equation},
\newblock Lett. Math. Phys. {\bf 11}, 247--245 (1986).

\bibitem{Fad.Res.Tak:89}
L.~D. Faddeev, N.~Yu. Reshetikhin, and L.~Takhtajan,
\newblock {\em Quantization of Lie groups and Lie algebras},
\newblock Alg. i Anal. {\bf 1}, 178--206 (1989),
\newblock (Leningrad Math. J. {\bf 1}, 193--225 (1990)).

\bibitem{Cel.Gia.Sor.Tar:91}
E.~Celeghini, R.~Giachetti, E.~Sorace, and M.~Tarlini,
\newblock {\em The quantum Heisenberg group $H(1)_q$; The three dimensional
  quantum group $E(3)_q$ and its $R$ matrix},
\newblock J. Math. Phys. {\bf 32}, 1155--1158, 1159--1165 (1991).

\bibitem{Cel.Gia.Sor.Tar:92}
E.~Celeghini, R.~Giachetti, E.~Sorace, and M.~Tarlini,
\newblock {\em Contractions of quantum groups}, vol. 1510 of {\em Lec. Notes
  Math.}, page 221,
\newblock Springer-Verlag, 1992.

\bibitem{Ino.Wig:53}
E.~\.In\"on\"u and E.~P. Wigner,
\newblock {\em On the contractions of groups and their representations},
\newblock Prot. Nat. Acad. Sci {\bf 39}, 510 (1953).

\bibitem{Luk.Now.Rue.Tol:91}
J.~Lukierski, A.~Nowicki, H.~Ruegg, and V.~N. Tolstoy,
\newblock {\em $q$ deformation of Poincar\'e algebra},
\newblock Phys. Lett. {\bf B264}, 331--338 (1991).

\bibitem{Luk.Rue.Tol:94}
J.~Lukierski, H.~Ruegg, and V.~N. Tolstoy,
\newblock {\em $\kappa$-Quantum Poincar\'{e} 1994},
\newblock in J.~Lukierski, Z.~Popowicz, and J.~Sobczyk, editors, {\em Quantum
  groups: formalism and applications}, page 359, PWN, 1994.

\bibitem{Mon.Pat:91}
M.~de~Montigny and J.~Patera,
\newblock {\em Discrete and continuous graded contractions of Lie algebras and
  superalgebras},
\newblock J. Phys. {\bf A24}, 525--547 (1991).

\bibitem{Moo.Pat:91}
R.~V. Moody and J.~Patera,
\newblock {\em Discrete and continuous graded contractions of representations
  of Lie algebras},
\newblock J. Phys. {\bf A24}, 2227--2257 (1991).

\bibitem{Bal.Her.Olm.San:94}
A.~Ballesteros, F.~J. Herranz, M.~A. del Olmo, and M.~Santander,
\newblock {\em $4-D$ quantum affine algebras and space-time $q$ symmetries},
\newblock J. Math. Phys. {\bf 35}, 4928--4940 (1994).

\bibitem{Bal.Her.Olm.San:95}
A.~Ballesteros, F.~J. Herranz, M.~A. del Olmo, and M.~Santander,
\newblock {\em Quantum algebras for maximal motion groups of $n$-dimensional
  flat spaces},
\newblock Lett. Math. Phys. {\bf 33}, 273--281 (1995).

\bibitem{Som:10}
D.~M.~Y. Sommerville,
\newblock Proc. Edinburgh Math. Soc. {\bf 28}, 25 (1910-11).

\bibitem{Yag.Roz.Yas:66}
I.~M. Yaglom, B.~A. Rozenfel'd, and E.~U. Yasinskaya,
\newblock Sov. Math. Surveys {\bf 19}, 49 (1966).

\bibitem{Maj:90}
S.~Majid,
\newblock {\em Physics for algebraists: Non-commutative and non-cocommutative
  Hopf algebras by a bicrossproduct construction},
\newblock J. Algebra {\bf 130}, 17--64 (1990).

\bibitem{Maj:90b}
S.~Majid,
\newblock {\em More examples of bicrossproduct and double cross product Hopf
  algebras},
\newblock Isr. J. Math. {\bf 72}, 133--148 (1990).

\bibitem{Maj:95}
S.~Majid,
\newblock {\em Foundations of quantum group theory},
\newblock Camb. Univ. Press, 1995.

\bibitem{Azc.Bue:96}
J.~A. de~Azc\'arraga and J.~C. {P\'erez Bueno},
\newblock {\em Contractions, Hopf algebra extensions and covariant differential
  calculus},
\newblock in B.~Jancewicz and J.~Sobczyk, editors, {\em From field theory to
  quantum groups}, pages 3--27, World Scientific, 1996.

\bibitem{Maj.Rue:94}
S.~Majid and H.~Ruegg,
\newblock {\em Bicrossproduct structure of $\kappa$-Poicar\'e group and
  non-commutative geometry},
\newblock Phys. Lett. {\bf B334}, 348--354 (1994).

\bibitem{Azc.Olm.Bue.San:97}
J.~A. de~Azc\'arraga, M.~A. del Olmo, J.~C. {P\'erez Bueno}, and M.~Santander,
\newblock {\em Graded contractions and bicrossproduct structure of deformed
  inhomogeneous algebras},
\newblock J. Phys. {\bf A30}, 3069--3086 (1997);

\bibitem{Per:97}
J.~C. {P\'erez Bueno},
\newblock {\em Bicrossproduct structure and graded contractions of deformed
  algebras},
\newblock Czech. J. Phys. {\bf 47}, 1275--1282 (1997).

\bibitem{Her.San:96b}
F.~J. Herranz and M.~Santander,
\newblock {\em The general solution of the real $Z_2^{\otimes N}$ graded
  contractions of $so(N+1)$},
\newblock J. Phys. {\bf A29}, 6643--6652 (1996).

\bibitem{Azc.Bue:96b}
J.~A. de~Azc\'arraga and J.~C. {P\'erez Bueno},
\newblock {\em Deformed and extended Galilei group Hopf algebras},
\newblock J. Phys. {\bf A29}, 6353--6362 (1996).

\end{thebibliography}

\end{document}